\documentclass[onecolumn]{autart}    

\usepackage{graphicx}
\usepackage{amsmath}
\usepackage{amsfonts}
\usepackage{amssymb}
\usepackage{etoolbox}
\usepackage{bm}
\usepackage{float}
\usepackage{epsfig}
\usepackage{multirow,bigdelim}
\usepackage{scalerel}
\usepackage{cite}
\usepackage{color,soul}
\usepackage{enumitem}

\begin{document}

\begin{frontmatter}

\title{Convergent Economic Model Predictive Control through parameter-varying storage functions for dissipativity} 

\thanks[footnoteinfo]{This paper was not presented at any IFAC 
meeting. Corresponding author: Zihang Dong}

\author[Paestum]{Zihang Dong}\ead{zihang.dong14@imperial.ac.uk},    
\author[Paestum,Rome]{David Angeli}\ead{d.angeli@imperial.ac.uk},               
\author[Paestum]{Goran Strbac}\ead{g.strbac@imperial.ac.uk}  

\address[Paestum]{Department of Electrical and Electronic Engineering,\\ Imperial College London, UK}  
\address[Rome]{Dipartimento di Ingegneria dell'Informazione, Universit\`a di Firenze, Italy}             


\begin{abstract}
This paper presents a new concept of controlled dissipativity as an extension of the standard dissipativity property to systems with parameter-varying storage functions under the framework of economic model predictive control (EMPC). Based on this concept, an EMPC controller, integrated with the dissipation inequality constraints treating the storage function parameters as decision variables, is formulated and the associated recursive feasibility is ensured. Then, the asymptotic convergence to an optimal equilibrium in closed-loop, without requiring the standard dissipativity assumption, is enforced by trading it off with asymptotic performance. Under suitable conditions on the storage function, stability of EMPC with respect to the optimal equilibrium is achieved. The upper bound of asymptotic average closed-loop performance is also evaluated. Finally, illustrative examples by using the EMPC controller with terminal regional conditions are provided to show the effectiveness of our methods.
\end{abstract}

\end{frontmatter}

\section{Introduction}\label{sec:intro}
The purpose of most advanced control systems is to eﬃciently and reliably guide a process to a desired setpoint. Model predictive control (MPC), due to its effectiveness in handling large and multi-variable systems subject to constraints, has become a prominent tool in many applications. The optimal steady-state setpoint is usually determined by some other information management system that decides which steady-state is the most advantageous from the economic point of view. However, this hierarchical separation of system information and objectives may be no longer optimal or desirable for an increasing number of applications.

In recent years, an alternative approach, known as economic model predictive control (EMPC), has looked into the issue of directly addressing economic optimisation in real time, and to this end, adopts objective costs which are not required to be positive definite w.r.t the optimal equilibrium point \cite{angeli2011average, rawlings2012fundamentals, grune2013economic, ellis2014tutorial, faulwasser2018economic}. Under this scheme, various tools in the literature have been proposed (see \cite{amrit2011economic,heidarinejad2012economic,grune2014asymptotic,angeli2016theoretical} for deterministic setups and \cite{marquez2014min,lucia2014handling,bayer2014tube,dong2020homothetic} for robust settings). To derive closed-loop guarantees, references \cite{grune2013economic,muller2016economic} employ the turnpike property of the underlying optimal control problem aiming at establishing (practical) closed-loop convergence to the optimal steady-state or periodic orbit. If the optimal regime of operation can be determined a-priori, for example, by validating a (strict) dissipativity condition, then use of EMPC method with appropriate terminal constraints can operate the system optimally \cite{diehl2010lyapunov,zanon2016periodic,dong2018analysis}. Specifically, a dissipativity property is often assumed to allow redefining a rotated positive definite versions of the economic cost, and the stability analysis is then carried out by adopting terminal equality or inequality constraints as in standard tracking MPC \cite{angeli2011average, rawlings2012fundamentals}.

The classical notions of (strict) dissipativity property are proposed in \cite{angeli2011average,grune2013economic} to prove stability of EMPC in the case of an optimal steady-state. The relations between dissipativity, turnpike property and the optimality of steady-state operation have been discussed in several publications. More precisely, the sufficiency of dissipativity for optimal steady-state operation is discussed in \cite{angeli2011average}, and the necessity, under a mild additional controllability assumption, is proved in \cite{muller2014necessity}. The equivalence between dissipativity and turnpike property under suitable assumptions is elaborated in \cite{grune2016relation}. The relation between dissipativity property, turnpike property, and optimal stationary operation in a continuous-time setting is investigated in \cite{faulwasser2014turnpike}. In addition, recent publications have presented different notions of dissipativity for systems with optimal periodic solutions. In \cite{muller2016economic}, a $P$-step MPC scheme is adopted to establish near optimal performance of the closed-loop system. Reference \cite{zanon2016periodic} defines periodic dissipativity as the existence of $P$ storage functions fulfilling $P$ dissipation inequalities. The equivalences between these two characterizations of dissipativity w.r.t periodic orbits are established in \cite{kohler2018periodic}.

The above literature has shown that dissipativity is playing a central role in the design and closed-loop analysis of EMPC schemes. However, in general, verifying the required dissipativity property w.r.t an equilibrium is a challenging task, and there is no systematic technique for general nonlinear optimal control problems. Particularly, dissipativity w.r.t a periodic orbit or more complicated non-periodic solution rather than a stationary equilibrium is significantly more difficult to verify. In \cite{pirkelmann2019approximate}, sum-of-squares techniques are used to numerically approximate the storage functions and thereby validating the dissipativity. A computational approach by solving the optimisation problem with dissipation inequality constraints for the verification of dissipativity properties w.r.t. general sets is presented in \cite{berberich2020dissipativity}. Moreover, the conditions for strict dissipativity w.r.t steady-state considering convex linear quadratic (LQ) constrained optimal control problems are provided in \cite{grune2018turnpike}. The results of \cite{grune2018turnpike} are later on generalized to indefinite LQ-problems by employing the particular shape of quadratic constraints to characterize strict dissipativity w.r.t. steady-state \cite{berberich2018indefinite}. However, for general nonlinear systems and/or non-convex feasible regions and cost functions, solving the optimal periodic or more complex solution a-priori is challenging and the performance optimisation may result in non-converging behaviours, hence enforcing convergence to a steady-state by trading it off with asymptotic performance may be desirable \cite{angeli2011enforcing}. Furthermore, convergence to the best equilibrium may be required when the knowledge of whether dissipativity holds or not is not available, e.g., the cost function changes on-line due to the varying economic environment and, even if sometimes dissipativity might be fulfilled, this is not easy to be a-priori assessed. In this context, establishing an EMPC scheme with suitable terminal constraints and terminal cost w.r.t an equilibrium can be applied to operate the system in a suboptimal way, which we will address in the present paper. Specifically, the main contributions and novelties of this paper can be summarized as follows:
\begin{itemize}
    \item The standard dissipativity property usually employed to guarantee stability of EMPC is loosened to existence of a parameter dependent controlled storage function. 
    
    \item An EMPC controller with terminal region requirements, is developed based on the proposed controlled dissipation inequality constraint that treats the storage function parameters as decision variables. The associated recursive feasibility is guaranteed and theoretically proved.
    
    \item The closed-loop asymptotic behavior using the EMPC controller is analysed. The result shows that the system trajectory can converge to the optimal equilibrium under the EMPC scheme with controlled dissipativity instead of the standard dissipativity. Moreover, asymptotic stability of the optimal equilibrium is achieved by imposing suitable technical conditions on the storage function. Finally, the upper bound of asymptotic average cost is evaluated.
    
\end{itemize}

The remainder of this paper is organized as follows. Section \ref{sec:preliminary} introduces the basic notions and setup. Section \ref{sec:empc_equality} presents the formulation of EMPC optimisation problem based on the controlled dissipativity constraint. Meanwhile, the proof of recursive feasibility for the proposed control algorithm are demonstrated. The asymptotic stability and average performance of the closed-loop system adopting the EMPC controller are analyzed in Section \ref{sec:stability_performance}. Finally, illustrative examples are presented in Section \ref{sec:example} and Section \ref{sec:conclusion} concludes this paper.

\textit{Notations}: Let symbols $\mathbb{R}$ and $\mathbb{I}$ denote the sets of real numbers and integers, respectively. $\mathbb{R}_{\geq a}$ ($\mathbb{I}_{\geq a}$) denotes the real numbers (integers) greater or equal to $a$ and $\mathbb{I}_{[a,b]}$ denotes the integers $\{a,a+1,\cdots,b\}$ for some $a \in \mathbb{I}$ and $b \in \mathbb{I}_{\geq a}$. Given a vector $x$, $\|x\|$ denotes the 2-norm. The symbol $x_{a:b|t}$ is a sequence from sampling instants $a$ to $b$, considered at time $t$. A continuous function $\alpha: \mathbb{R}_{\geq 0} \to \mathbb{R}_{\geq 0}$ is of $class$ $\mathcal{K}$, if it is zero at zero and strictly increasing. It is of $class$ $\mathcal{K}_{\infty}$ if in addition, $\alpha(s) \to +\infty$ as $s \to +\infty$. A continuous function $\rho: \mathbb{R}^{n} \to \mathbb{R}$ is positive definite with respect to some point $z \in \mathbb{R}^{n}$ if $\rho(z) = 0$ and $\rho(x)>0$ for all $x \neq z$.

\section{Preliminaries}\label{sec:preliminary}

\subsection{Problem description}

We consider the finite-dimensional discrete time nonlinear system
\begin{equation}\label{nonlinearsystem}
    x_{t+1} = f(x_{t}, u_{t})
\end{equation}
with state $x \in \mathbb{X} \subseteq \mathbb{R}^{n}$, input $u \in \mathbb{U} \subseteq \mathbb{R}^{m}$ and continuous state transition map $f: \mathbb{X} \times\mathbb{U} \to \mathbb{X}$. The state and input of the system must fulfill the
pointwise-in-time state and input constraints
\begin{equation}\label{stateinputconstraints}
    (x_{t},u_{t}) \in \mathbb{Z} \qquad  \forall t \in \mathbb{I}_{\geq 0}
\end{equation}
for some compact set $\mathbb{Z} \subseteq \mathbb{X} \times \mathbb{U}$. Constraint \eqref{stateinputconstraints} implies the control constraint is possibly state-dependent, i.e.,
$u \in \mathbb{U}(x) = \{ u \in \mathbb{U} \,|\, (x,u) \in \mathbb{Z} \}$ and implicitly enforces the state constraint $x \in \{x \in \mathbb{X}\,|\, \mathbb{U}(x) \neq \emptyset\}$. If there are no coupled constraints, then the state-input constraints in \eqref{stateinputconstraints} become $x \in \mathbb{X}$ and $u \in \mathbb{U}$.

The objective we seek to optimise is the  following economic costs accumulated in the long term system operation
\begin{equation}\label{infiniteobjective}
    \sum^{\infty}_{t=0} \ell(x_{t},u_{t})
\end{equation}
subject to the dynamic constraints \eqref{nonlinearsystem} and state-input constraints \eqref{stateinputconstraints}. The function $\ell(x,u): \mathbb{Z} \to \mathbb{R}$, specifying the performance criterion to be minimized, is continuous and denotes the stage cost that may take arbitrary form coherently with the EMPC setup and need not be positive definite w.r.t any equilibrium state. With the notation adopted so far, the optimal feasible steady-state fulfils
\begin{equation}
    (x^{s},u^{s}) = \arg \min \{\ell(x,u)\,|\, (x,u)\in \mathbb{Z}, x=f(x,u)\}.
\end{equation}
and throughout this work we assume the solution of the steady-state problem is unique.

\subsection{Controlled dissipativity property}

Given the system dynamics \eqref{nonlinearsystem}, the constraint set $\mathbb{Z}$ and the stage cost $\ell(x,u)$, our aim is to investigate how to design a controller ensuring stable closed-loop system and showing convergence to an optimal equilibrium. The initial stability analysis of EMPC was presented in \cite{rawlings2008unreachable} for linear systems with strictly convex cost functionals. Subsequently, for a class of nonlinear systems and an arbitrary economic objective, reference \cite{diehl2010lyapunov}  defines the rotated cost-to-go as a candidate Lyapunov function under the assumption of strong duality. To relax the assumptions of \cite{diehl2010lyapunov}, a dissipativity notion in proving stability of EMPC is introduced in \cite{angeli2011average}. Let us first recall the definition of standard dissipativity.
\begin{defn}\textup{{\cite[Definition 4.1]{angeli2011average}}}\label{standarddissipativitydefinition}
A control system as in \eqref{nonlinearsystem} is dissipative w.r.t a supply rate $s: \mathbb{X}\times\mathbb{U} \to \mathbb{R}$ if there exists a function $\lambda: \mathbb{X} \to \mathbb{R}$ (the storage function) such that
\begin{equation}
    \lambda(f(x,u)) - \lambda(x) \leq s(x,u),
\end{equation}
for all $(x,u)\in \mathbb{Z} \subseteq \mathbb{X} \times \mathbb{U}$. If, in addition, a positive definite function $\rho: \mathbb{X} \to \mathbb{R}_{\geq 0}$ exists such that
\begin{equation}
    \lambda(f(x,u)) - \lambda(x) \leq s(x,u) - \rho(x-x^{s})
\end{equation}
then the system is said to be strictly dissipative.
\end{defn}

It is worth to remark that the supply rate introduced above is usually in the form of $s(x,u) = \ell(x,u) - \ell(x^{s},u^{s})$ where $\ell$ refers to the stage cost. The dissipativity property is of important because it establishes a relation between the economic cost and the system trajectories, and it allows to make qualitative statements about the optimal operation of a process on infinite horizons. In particular, \cite[Proposition 6.4]{angeli2011average} shows that dissipativity property is sufficient to imply that the system is optimally operated at steady-state, which can be stated as: for any solution such that $(x_{t},u_{t})\in \mathbb{Z}$ for all $t \in \mathbb{I}_{\geq 0}$, it holds
\begin{equation}\label{optimalsteadystateoperation}
    \liminf_{T \to \infty} \frac{\sum^{T-1}_{t=0} \ell(x_{t},u_{t})}{T} \geq \ell(x^{s},u^{s}).
\end{equation}
Inequality \eqref{optimalsteadystateoperation} indicates that no feasible state and input sequence pair can offer a better asymptotic average performance than the optimal steady-state cost.

Note that the above dissipativity condition could be regarded as a strong condition because the dissipation inequality must hold for all feasible state-input pairs and the storage function parameters are time-invariant. However, due to the plant’s nonlinearities and nonconvex cost functionals, the best operating regime for given system and constraints could actually fail to be an equilibrium, which means the dissipativity in Definition \ref{standarddissipativitydefinition} does not hold. 
Moreover, obtaining a periodic or even complex chaotic optimal solution is usually challenging and operating the system with a non-stationary control input could lead to unsafe behaviours, so one may trade off the economic performance by gradually steering the plant towards the best steady-state.
\begin{defn}\label{dissipativitydefinition}
The system \eqref{nonlinearsystem} is controlled dissipative with a parameter-varying storage function $\lambda: \Theta \times \mathbb{X} \to \mathbb{R}$ if for all parameter $\theta \in \Theta$ and all state $x \in \mathbb{X}$ of the parameter varying storage function $\lambda(\theta,x)$, there exists $\theta^{+} \in \Theta$ and $u \in \mathbb{U}$ such that the following inequality holds,
\begin{equation}\label{controlleddissipationinequality}
    \lambda(\theta^{+},f(x,u)) - \lambda(\theta,x) \leq s(x,u),
\end{equation}
where the supply rate $s(x,u) = \ell(x,u) - \ell(x^{s},u^{s})$. If, in addition, there exists a positive definite function $\rho: \mathbb{X} \to \mathbb{R}_{\geq 0}$ such that
\begin{equation}\label{dissipationinequality0}
    \lambda(\theta^{+},f(x,u)) - \lambda(\theta,x) \leq \ell(x,u) - \ell(x^{s},u^{s}) - \rho( x-x^{s} )
\end{equation}
then the system is strictly controlled dissipative.
\end{defn}

\begin{rem}
We remark that the controlled dissipativity in Definition \ref{dissipativitydefinition} was first proposed in \cite{dong2018analysis} for non parametric storage functions. 
In \cite{dong2018analysis}, however, controlled dissipativity property is only adopted for the terminal penalty function
whereas controlled dissipativity in Definition \ref{dissipativitydefinition} is enforced along open-loop solutions along the lines of the standard dissipativity. 
Particularly, the controlled dissipativity with parameter-varying storage functions in \eqref{controlleddissipationinequality} and \eqref{dissipationinequality0} is a relaxed condition compared to the standard dissipativity in Definition \ref{standarddissipativitydefinition} which allows to retain stability of closed-loop EMPC
even when the considered dynamics and stage-cost are not dissipative. The parameters of the storage function are allowed to vary in time and will be considered as decision variables in the EMPC algorithm design.
\end{rem}

Then, we make the following assumptions.
\begin{assum}\label{assum:compactTheta}
The set $\Theta$ of storage function parameters is compact.
\end{assum}

\begin{assum}\label{assum:continuousstoragefunction}
The storage function $\lambda(\theta,x)$ is continuous on $\Theta\times \mathbb{X}$.
\end{assum}

Since the stage cost $\ell(x,u)$ is continuous and the set $\mathbb{Z}$ is compact, the supply rates on the right hand side of \eqref{controlleddissipationinequality} and \eqref{dissipationinequality0} are lower bounded. Additionally, based on Assumption \ref{assum:compactTheta} and \ref{assum:continuousstoragefunction}, the storage function $\lambda(\theta,x)$ is also lower bounded and can be shifted by constant terms. Hence, the larger set $\Theta$ implies that controlled dissipativity is less restrictive and the system has more ``spare time" to seek for improved economic performance during the transient period (see examples in Section \ref{sec:example}).


\section{Economic MPC with terminal constraints}\label{sec:empc_equality}

In this section, an EMPC design with terminal constraint and penalty for the controlled dissipative system is presented. We firstly provide the detailed formulation and description of the optimisation problem. Then, the associated recursive feasibility is proved based on a constructed suboptimal solution.

\subsection{Algorithm formulation}\label{empc:terminalequality}

In general, the standard dissipativity property in Definition \ref{standarddissipativitydefinition} is assumed to develop an EMPC controller that is able to guarantee asymptotic stability for the optimal steady-state. When a steady-state is not the optimal regime of operation, the standard dissipativity does not hold and asymptotic stability of closed-loop EMPC w.r.t the steady-state is lost. In this context, the following assumption of controlled dissipativity is made
\begin{assum}\label{assumptioncontrolleddissipative}
System \eqref{nonlinearsystem} is strictly controlled dissipative with respect to the supply rate $s(x,u) = \ell(x,u)-\ell(x^{s},u^{s})$.
\end{assum}

Next, to avoid the nonconvex infinite dimensional problem in \eqref{infiniteobjective}, we here define an optimisation problem over a sufficiently long, but finite horizon, i.e., the EMPC optimisation problem over a finite prediction horizon $N \in \mathbb{I}_{[1,+\infty)}$ at any sampling time instant $t$ is formulated as
\begin{subequations}\label{Opt:main}
\begin{align}
    \min_{\bm{\theta}, \bm{x}, \bm{u}} &\,\, J_{N}(x_{t:t+N|t},u_{t:t+N-1|t}) \tag{\ref{Opt:main}}\\
    s.t.\quad & x_{t|t} = x_{t},  \label{Opt:initial}\\
    & \theta_{t|t} = \theta_{t},  \label{Opt:initialtheta}\\
    & x_{k+1|t} = f(x_{k|t}, u_{k|t}), \,\forall k \in \mathbb{I}_{[t,t+N-1]}  \label{Opt:dynamic}\\
    & (x_{k|t}, u_{k|t}) \in \mathbb{Z},\,\forall k \in \mathbb{I}_{[t,t+N-1]} \label{Opt:feasible}\\
    & \theta_{k|t} \in \Theta,\,\forall k \in \mathbb{I}_{[t,t+N]}  \label{Opt:feasibletheta}\\
    & \lambda(\theta_{k+1|t}, x_{k+1|t}) - \lambda(\theta_{k|t}, x_{k|t}) + \rho(x_{k|t}-x^{s})  \nonumber \\
    &\qquad \qquad \leq \ell(x_{k|t},u_{k|t}) - \ell(x^{s},u^{s}),\,\forall k \in \mathbb{I}_{[t,t+N-1]}  \label{Opt:dissipativity}\\
    & x_{t+N|t} \in \mathbb{X}_{f} \label{Opt:terminal}
\end{align}
\end{subequations}
where $u_{t:t+N-1|t} = \{u_{t|t}, u_{t+1|t}, \cdots u_{t+N-1|t}\}$ is the control sequence, $x_{t:t+N|t} = \{x_{t|t}, x_{t+1|t}, \cdots x_{t+N-1|t}, x_{t+N|t}\}$ is the associated state trajectory, and $\theta_{t:t+N|t} = \{\theta_{t|t}, \theta_{t+1|t}, \cdots \theta_{t+N|t}\}$ is the sequence of storage function parameters within the prediction horizon.
The objective function in \eqref{Opt:main} is explicitly defined as
\begin{equation}\label{objective}
    J_{N}(x_{t:t+N|t},u_{t:t+N-1|t}) = \sum^{t+N-1}_{k=t} \ell(x_{k|t},u_{k|t}) + V_{f}(x_{t+N|t}).
\end{equation}
with the continuous terminal penalty function $V_{f}(\cdot)$ defined on the terminal region $\mathbb{X}_{f}$ as shown in \eqref{Opt:terminal}.

Constraints \eqref{Opt:initial} and \eqref{Opt:initialtheta} are the initial conditions, in particular, \eqref{Opt:initialtheta} requires the parameters of storage function at current time $t$ be fixed at the optimal predicted value $\theta_{t} = \theta^{*}_{t|t-1}$ of time $t-1$.
This constraint, while not needed for recursive feasibility, will be used to show the equivalence between MPC formulations using rotated cost functions and the original one (as customary for stability analysis in EMPC).
Conditions \eqref{Opt:dynamic} represents the nonlinear system dynamics. The feasibility of system state-input pair and storage function parameters are ensured by \eqref{Opt:feasible} and \eqref{Opt:feasibletheta}, respectively. The dissipation inequality is imposed by the constraint \eqref{Opt:dissipativity}. Thanks to Assumption \ref{assumptioncontrolleddissipative}, it is feasible to obtain solutions of control inputs and storage function parameters over the prediction horizon. In particular, a positive definite term w.r.t the distance $x_{k|t}-x^{s}$ is used to enforce the system is strict dissipative and it also determines the convergence speed which will be assessed in Section \ref{sec:example}. Finally, \eqref{Opt:terminal} represents the terminal region requirement.

\begin{rem}
Although the constraints \eqref{Opt:dissipativity} are included in the optimisation problem \eqref{Opt:main}, it does not impose additional restrictions compared to the standard EMPC formulation if dissipativity is fulfilled with a storage function of the form $\lambda(\theta,x) $ for some constant value of $\theta$. Notice that in the standard EMPC setup a dissipativity condition w.r.t all $(x,u)\in \mathbb{Z}$ is usually assumed in the stability analysis. In our proposed method, the strong requirement of dissipativity condition is relaxed to satisfaction of dissipation inequalities along a certain state and input trajectory determined by the optimisation problem.
\end{rem}

To establish recursive feasibility and analyse system performance, the following assumption regarding terminal region $\mathbb{X}_{f}$ and penalty function $V_{f}$ is made:
\begin{assum}\label{assum:terminal}
There exists a compact terminal region $\mathbb{X}_{f} \subseteq \mathbb{X}$ with $x^{s}$ as an interior point, a terminal policy $\kappa_{f}(x): \mathbb{X} \to \mathbb{U}$ with $\kappa_{f}(x^{s}) = u^{s}$, and a continuous terminal penalty function $V_{f}(\cdot): \mathbb{X}_{f} \to \mathbb{R}$, such that for all $x \in \mathbb{X}_{f}$ the following hold:\\
(i)\,\, $(x, \kappa_{f}(x)) \in \mathbb{Z}$;\\
(ii)\, $f(x, \kappa_{f}(x)) \in \mathbb{X}_{f}$;\\
(iii) $V_{f}(f(x,\kappa_{f}(x))) - V_{f}(x) \leq \ell(x^{s},u^{s}) - \ell(x,\kappa_{f}(x))$.
\end{assum}

Denote the optimal solution of problem \eqref{Opt:main} at time $t \in \mathbb{I}_{\geq 1}$ as
\begin{subequations}\label{optimalsolution}
\begin{equation}
    u^{*}_{t:t+N-1|t} = \{u^{*}_{t|t}, u^{*}_{t+1|t}, \cdots, u^{*}_{t+N-1|t}\}
\end{equation}
\begin{equation}
    x^{*}_{t:t+N|t} = \{x^{*}_{t|t}, x^{*}_{t+1|t}, \cdots, x^{*}_{t+N-1|t}, x^{*}_{t+N|t}\}
\end{equation}
\begin{equation}
    \theta^{*}_{t:t+N|t} = \{\theta^{*}_{t|t}, \theta^{*}_{t+1|t}, \cdots, \theta^{*}_{t+N-1|t}, \theta^{*}_{t+N|t}\}.
\end{equation}
\end{subequations}
These sequences are computed from a warm-start that is a feasible and suboptimal solution developed at the previous optimisation step, viz.
\begin{subequations}\label{initialsequence}
\begin{equation}
    \tilde{u}_{t:t+N-1|t} = \{u^{*}_{t|t-1}, \cdots, u^{*}_{t+N-2|t-1}, \kappa_{f}(x^{*}_{t+N-1|t-1})\}
\end{equation}
\begin{equation}
    \tilde{x}_{t:t+N|t} = \{x^{*}_{t|t-1}, \cdots, x^{*}_{t+N-1|t-1}, x^{*}_{t+N|t-1}\}
\end{equation}
\begin{equation}
    \tilde{\theta}_{t:t+N|t} = \{\theta^{*}_{t|t-1}, \cdots, \theta^{*}_{t+N-1|t-1}, \theta^{*}_{t+N|t-1}\}
\end{equation}
\end{subequations}
where $x^{*}_{t+N|t-1} = f(x^{*}_{t+N-1|t-1},\kappa_{f}(x^{*}_{t+N-1|t-1}))$.


With the optimal solution of \eqref{optimalsolution}, the optimal objective function can be obtained as
\begin{equation}
    V_{N}(\theta_{t},x_{t}) = J_{N}(x^{*}_{t:t+N|t},u^{*}_{t:t+N-1|t}),
\end{equation}
Following the manner of feedback algorithm in EMPC, the optimal control implemented to system \eqref{nonlinearsystem} is $u^{*}_{t|t}$ and the closed-loop system dynamic is
\begin{equation}\label{closedloopsystem}
    x_{t+1} = f(x_{t},u^{*}_{t|t}).
\end{equation}
In addition, the optimal value $\theta^{*}_{t+1|t}$ will be adopted as the initial condition \eqref{Opt:initialtheta} of the optimisation problem at the next sampling time instant to ensure recursive feasibility. 

\subsection{Recursive feasibility}

We now turn our attention to the recursive feasibility of EMPC optimisation problem \eqref{Opt:main}. For an initial real system state $x$ and an initial storage function parameter $\theta$, let us denote the admissible set of solutions of \eqref{Opt:main} fulfilling \eqref{Opt:initial}-\eqref{Opt:terminal} by $\mathcal{Z}_{N}(x,\theta)$, then the set of admissible states and parameters is defined as
\begin{equation}\label{regionofattraction}
    \mathcal{X}_{N} = \{(x,\theta) \in \mathbb{X} \times \Theta\,|\, \exists (\bm{x}, \bm{u}, \bm{\theta}) \in \mathcal{Z}_{N}(x,\theta)\}.
\end{equation}

While the terminal equality condition is used, i.e., $x_{t+N|t} = x^{s}$, the suboptimal terminal storage function parameter can be chosen as $\theta^{*}_{t+N|t-1} = \theta^{*}_{t+N-1|t-1}$ such that the resulting controlled dissipativity inequality \eqref{Opt:dissipativity} is trivially satisfied for $u_{t+N-1|t}=u^s$. However, the terminal region constraint \eqref{Opt:terminal} only enforces $x_{t+N|t} \in \mathbb{X}_{f}$ and thereby the previous considerations do not apply, so that
choosing $\theta^{*}_{t+N|t-1} = \theta^{*}_{t+N-1|t-1}$ is not always feasible. In this case, the following assumption concerning the terminal storage function parameter is required. 
\begin{assum}\label{assumptiondissipative}
For all state $x \in \mathbb{X}_{f}$ and all $\theta \in \Theta$, there exist parameter $\theta^{+}\in \Theta$ and positive definite $\rho(\cdot)$ such that the following strict dissipation inequality under terminal control policy $\kappa_{f}(x)$ holds,
\begin{equation}
    \lambda(\theta^{+},f(x,\kappa_{f}(x))) - \lambda(\theta,x)  \leq \ell(x,\kappa_{f}(x)) - \ell(x^{s},u^{s}) - \rho(x-x^{s})
\end{equation}
\end{assum} 

Now, we are ready to claim the recursive feasibility
of the optimisation problem \eqref{Opt:main}.
\begin{prop}\label{recursivefeasibilityequality}
Let Assumption \ref{assum:compactTheta}, \ref{assum:continuousstoragefunction}, \ref{assum:terminal}(i)(ii) and \ref{assumptiondissipative} hold. Then, for any initial state and storage function parameter $(x,\theta) \in \mathcal{X}_{N}$, the economic MPC optimisation problem \eqref{Opt:main} is recursively feasible.
\end{prop}
\begin{pf}
The main idea to prove recursive feasibility is to explicitly construct feasible solutions for optimisation problem \eqref{Opt:main} at the current time instant $t \in \mathbb{I}_{\geq 1}$, given the feasible and optimal solution generated at time $t-1$. Specifically, the sequences in \eqref{initialsequence} will be shown as a feasible solution of problem \eqref{Opt:main}.

Since the sequence \eqref{initialsequence} are constructed from the optimisation problem at time $t-1$, they fulfill the initial conditions \eqref{Opt:initial}-\eqref{Opt:initialtheta}, the system dynamic equality constraint \eqref{Opt:dynamic} for $k\in \mathbb{I}_{t:t+N-1}$, the state-input and storage function parameter constraint \eqref{Opt:feasible}-\eqref{Opt:feasibletheta} for $k\in \mathbb{I}_{t:t+N-1}$, and the dissipation inequality for $k\in \mathbb{I}_{t:t+N-2}$, automatically. Regarding the dissipation inequality at $k=t+N-1$, based on Assumption \ref{assumptiondissipative}, there always exist storage function parameter $\theta^{*}_{t+N|t-1} \in \Theta$ such that the controlled dissipativity condition is satisfied under the terminal control input $\kappa_{f}(x^{*}_{t+N-1|t-1})$, therefore, we can choose $\theta^{*}_{t+N|t-1} = \theta^{*}_{t+N-1|t-1}$ such that the dissipation inequalities \eqref{Opt:dissipativity} are also satisfied. Furthermore, By imposing the terminal constraint \eqref{Opt:terminal} at time $t-1$, it holds $x^{*}_{t+N-1|t-1} \in \mathbb{X}_{f}$ because of the terminal conditions in Assumption \ref{assum:terminal}.

Therefore, the optimisation problem at current time $t$ has at least one solution, and recursive feasibility of this optimisation problem is ensured. \qed
\end{pf}

\section{Closed-loop analysis and performance}\label{sec:stability_performance}

This section explores the asymptotic stability and average economic performance of the closed-loop system under EMPC control actions. Based on the proposed EMPC controller, the first part is to prove the asymptotic convergence and stability of the optimal steady-state, and then an upper bound of the asymptotic average performance is derived.

\subsection{Convergence and stability analysis}

As discussed in Section \ref{sec:intro}, dissipativity is sufficient for ensuring asymptotic stability of the optimal steady-state. If the standard dissipativity condition is not satisfied, convergence can be enforced by imposing average constraints. Specifically, for an economic setting with average constraints, a first preliminary result of closed-loop convergence was obtained by \cite[Remark 6.5]{angeli2011average} using the terminal equality constraint. Reference \cite{muller2014convergence} extended the convergence result under average constraints to an EMPC scheme with more general terminal region and terminal penalty function. However, the approach through average constraints, does not allow to warrant the asymptotic stability of $x^{s}$. 
To close this gap, we prove that the use of controlled dissipativity constraint can enforce convergence and asymptotic stability by allowing to trade-off convergence to optimal equilibria and economic performance, while retaining Lyapunov stability.
\begin{thm}\label{thm:convergence}
Suppose Assumption \ref{assum:compactTheta}, \ref{assum:continuousstoragefunction}, \ref{assumptioncontrolleddissipative}, \ref{assum:terminal}, and \ref{assumptiondissipative} are satisfied. Then the solution of EMPC closed-loop system \eqref{closedloopsystem} with initial condition $(x_{0},\theta_{0}) \in \mathcal{X}_{N}$ asymptotically converges to $x^{s}$, i.e., $\lim_{t\to\infty} x_{t} = x^{s}$. 
\end{thm}
\begin{pf}
We first introduce the rotated stage
cost function and rotated terminal cost function as
\begin{equation*}
\begin{aligned}
    L(\theta,\theta^{+},x,u) &= \ell(x,u) + \lambda(\theta,x) - \lambda(\theta^{+},f(x,u)) - \ell(x^{s}, u^{s})\\
    \bar{V}_{f}(\theta,x) &= V_{f}(x) + \lambda(\theta,x)
\end{aligned}
\end{equation*}
Then, according to the controlled dissipativity in Assumption \ref{assumptioncontrolleddissipative}, it implies that for all $x\in \mathbb{X}$ and all $\theta \in \Theta$ there exists $u \in \mathbb{U}$ and $\theta^{+} \in \Theta$ such that the rotated stage cost $L(\theta,\theta^{+},x,u)\geq \rho(x-x^{s})$. Thus, it is lower bounded by a \textit{class} $\mathcal{K}$ function $\gamma(\cdot)$, i.e., $ L(\theta,\theta^{+},x,u) \geq \rho(x-x^{s}) \geq \gamma(\|x-x^{s}\|)$. For all terminal state $x \in \mathbb{X}_{f}$, the rotated terminal cost function fulfils
\begin{equation}\label{rotatedterminalinequality}
\begin{aligned}
    \bar{V}_{f}(\theta^{+},x^{+}) - \bar{V}_{f}(\theta,x) =& V_{f}(x^{+}) + \lambda(\theta^{+},x^{+}) - V_{f}(x) - \lambda(\theta,x)\\
    \leq & \ell(x^{s},u^{s}) - \ell(x,\kappa_{f}(x)) + \lambda(\theta^{+},x^{+}) - \lambda(\theta,x)\\
    =& -L(\theta,\theta^{+},x,\kappa_{f}(x))
\end{aligned}
\end{equation}
where $x^{+}=f(x,\kappa_{f}(x))$ and the inequality is derived according to Assumption \ref{assum:terminal}.

Next, for initial condition $(x_{0},\theta_{0}) \in \mathcal{X}_{N}$ and terminal state $x_{t+N|t} \in \mathbb{X}_{f}$, the rotated objective function at time $t$ is expressed as
\begin{equation*}
\begin{aligned}
    \bar{J}_{N}&(\theta_{t:t+N|t},x_{t:t+N|t},u_{t:t+N-1|t})\\
    =& \sum^{t+N-1}_{k=t} L(\theta_{k|t},\theta_{k+1|t},x_{k|t},u_{k|t}) + \bar{V}_{f}(\theta_{t+N|t},x_{t+N|t})\\
    =& \bigg[\ell(x_{t|t},u_{t|t})+\lambda(\theta_{t|t},x_{t|t})-\lambda(\theta_{t+1|t},x_{t+1|t})- \ell(x^{s}, u^{s})\bigg] + \bigg[\ell(x_{t+1|t},u_{t+1|t})+\lambda(\theta_{t+1|t},x_{t+1|t})-\lambda(\theta_{t+2|t},x_{t+2|t})\\ 
    & - \ell(x^{s}, u^{s})\bigg] +  \cdots +  \bigg[\ell(x_{t+N-1|t},u_{t+N-1|t}) +\lambda(\theta_{t+N-1|t},x_{t+N-1|t})-\lambda(\theta_{t+N|t},x_{t+N|t}) - \ell(x^{s}, u^{s})\bigg]\\& +  \bigg[ V_{f}(x_{t+N|t}) + \lambda(\theta_{t+N|t},x_{t+N|t}) \bigg]\\
    =& J_{N}(x_{t:t+N|t},u_{t:t+N-1|t})  - N \ell(x^{s},u^{s}) + \lambda(\theta_{t|t},x_{t|t})
\end{aligned}
\end{equation*}
Since $\lambda(\theta_{t|t},x_{t|t})$ is evaluated at the given initial condition \eqref{Opt:initial}-\eqref{Opt:initialtheta} and $N \ell(x^{s},u^{s})$ is independent of the decision variables, the two objective functions $\bar{J}_{N}$ and $J_{N}$ differ by a constant, therefore the rotated cost function $\bar{J}_{N}$ subject to constraints \eqref{Opt:initial}-\eqref{Opt:terminal} admits the same optimal solution as the original optimisation problem \eqref{Opt:main}. Based on the assumptions of compactness of sets and continuity of functions, solutions exist for both \eqref{Opt:main} and its rotated optimisation problem.

Next, based on the optimal rotated value function, i.e., $\bar{V}_{N}(\theta_{t},x_{t})=\bar{J}_{N}(\theta^{*}_{t:t+N|t},x^{*}_{t:t+N|t},u^{*}_{t:t+N-1|t})$, the following inequality can be obtained
\begin{equation*}
\begin{aligned}
    \bar{J}_{N}&(\tilde{\theta}_{t:t+N|t},\tilde{x}_{t:t+N|t},\tilde{u}_{t:t+N-1|t})\\
    =& \bar{V}_{N}(\theta_{t-1},x_{t-1}) - L(\theta^{*}_{t-1|t-1},\theta^{*}_{t|t-1},x^{*}_{t-1|t-1},u^{*}_{t-1|t-1}) - \bar{V}_{f}(\theta^{*}_{t+N-1|t-1},x^{*}_{t+N-1|t-1})\\
    &+L(\theta^{*}_{t+N-1|t-1},\theta^{*}_{t+N|t-1},x^{*}_{t+N-1|t-1},\kappa_{f}(x^{*}_{t+N-1|t-1})) + \bar{V}_{f}(\theta^{*}_{t+N|t-1},f(x^{*}_{t+N-1|t-1},\kappa_{f}(x^{*}_{t+N-1|t-1})))\\
    \leq& \bar{V}_{N}(\theta_{t-1},x_{t-1}) - L(\theta^{*}_{t-1|t-1},\theta^{*}_{t|t-1},x^{*}_{t-1|t-1},u^{*}_{t-1|t-1}).
\end{aligned}
\end{equation*}
where the last inequality holds because of \eqref{rotatedterminalinequality}.

Since $\bar{V}_{N}(\theta_{t},x_{t}) \leq \bar{J}_{N}(\tilde{\theta}_{t:t+N|t},\tilde{x}_{t:t+N|t},\tilde{u}_{t:t+N-1|t})$, it follows that
\begin{equation*}
    \bar{V}_{N}(\theta_{t},x_{t}) - \bar{V}_{N}(\theta_{t-1},x_{t-1}) \leq -L(\theta^{*}_{t-1|t-1},\theta^{*}_{t|t-1},x^{*}_{t-1|t-1},u^{*}_{t-1|t-1}) \leq -\rho(x-x^{s}) \leq -\gamma(\|x_{t}-x^{s}\|),
\end{equation*}
thus the optimal rotated value function $\bar{V}_{N}(\theta_{t},x_{t})$ is non-increasing along the trajectory of the closed-loop system \eqref{closedloopsystem}. Moreover, due to the continuity of functions and the compactness of sets, the sequence $\bar{V}_{N}$ is bounded from below and $\bar{V}_{N}(x_{0},\theta_{0})$ is finite. Hence, the positive definite functions $\sum^{\infty}_{t=0} \rho(x_{t}-x^{s})$ converges which implies that $x_{t}$ asymptotically converges, i.e., $\lim_{t\to\infty} x_{t} = x^{s}$. \qed
\end{pf}

Although Theorem \ref{thm:convergence} shows asymptotic convergence to the optimal steady-state $x^{s}$, this is not necessarily asymptotically stable because Lyapunov stability property has not been guaranteed. In particular, a uniform upper bound on the rotated value function $\bar{V}_{N}$ needs to be derived for all $\theta \in \Theta$. Note that using the rotated value function as a Lyapunov function requires $\bar{V}_{f}(\theta,x^{s}) = 0, \forall \theta \in \Theta$. Accordingly, the following assumption is made:
\begin{assum}\label{assum:stability}
The storage function satisfies $\lambda(\theta,x^{s}) = 0, \forall \theta \in \Theta$.
\end{assum}

We can now claim the main result on asymptotic stability of closed-loop EMPC as below
\begin{cor}\label{theoreminequality}
Let Assumption \ref{assum:compactTheta}, \ref{assum:continuousstoragefunction}, \ref{assumptioncontrolleddissipative}, \ref{assum:terminal}, \ref{assumptiondissipative}, and \ref{assum:stability} hold, the optimal steady-state $(x^{s}, u^{s})$ is asymptotically stable with a region of attraction $\mathcal{X}_{N}$ by using the EMPC controller \eqref{Opt:main}.
\end{cor}

\begin{pf}
In addition to the non-increasing property of $\bar{V}_{N}(\theta_{t},x_{t})$ analysed in the proof of Theorem \ref{thm:convergence}, we need to provide lower and upper bounds for the function $\bar{V}_{N}(\theta_{t},x_{t})$ to use it as a Lyapunov function. Because the rotated stage cost is non-negative and the rotated terminal cost can be shifted to positive by a constant value, the optimal rotated objective function fulfills
\begin{equation*}
    \bar{V}_{N}(\theta_{t},x_{t}) \geq L(\theta^{*}_{t|t},\theta^{*}_{t+1|t},x_{t},u^{*}_{t|t}) \geq \gamma(\|x_{t}-x^{s}\|).
\end{equation*}
On the other hand, similar to \cite[Proposition 2.18]{rawlings2017model}, from the dynamic programming recursion, $\forall x_{t} \in \mathbb{X}_{f}$, it holds 
\begin{equation*}
    \bar{V}_{1}(\theta_{t},x_{t}) \leq L(\theta^{*}_{t|t},\theta^{*}_{t+1|t},x_{t},\kappa_{f}(x_{t})) + \bar{V}_{f}(\theta^{*}_{t+1|t},f(x_{t},\kappa_{f}(x_{t})))
    \leq \bar{V}_{f}(\theta^{*}_{t|t},x_{t})
\end{equation*}
where the first inequality holds because of the use of suboptimal control policy $\kappa_{f}(x_{t})$ and the second is the rotated terminal inequality. Then, by induction, for fixed prediction horizon $N$ and $x_{t} \in \mathbb{X}_{f}$, it yields
\begin{equation*}
\begin{aligned}
    \bar{V}_{N}(\theta_{t},x_{t}) \leq& L(\theta^{*}_{t|t},\theta^{*}_{t+1|t},x_{t},\kappa_{f}(x_{t})) + \bar{V}_{N-1}(\theta^{*}_{t+1|t},f(x_{t},\kappa_{f}(x_{t})))\\
    \leq& L(\theta^{*}_{t|t},\theta^{*}_{t+1|t},x_{t},\kappa_{f}(x_{t})) + \bar{V}_{f}(\theta^{*}_{t+1|t},f(x_{t},\kappa_{f}(x_{t})))\\
    \leq& \bar{V}_{f}(\theta^{*}_{t|t},x_{t}).
\end{aligned}
\end{equation*}
Moreover, based on Assumption \ref{assum:stability} and the uniform continuity property of functions on compact sets, we can obtain $\bar{V}_{f}(\theta^{*}_{t|t},x_{t}) \leq \gamma_{f}(\|x_{t}-x^{s}\|), \forall x_{t} \in \mathbb{X}_{f}$, where $\gamma_{f}(\cdot)$ is a \textit{class} $\mathcal{K}_{\infty}$ function. According to \cite[Proposition 2.16]{rawlings2017model}, the local upper bound on $\mathbb{X}_{f}$ implies the existence of upper bound of $\bar{V}_{N}(\theta_{t},x_{t})$ to $\mathcal{X}_{N}$, i.e., $\bar{V}_{N}(\theta_{t},x_{t}) \leq \hat{\gamma}(\|x_{t}-x^{s}\|), \forall (x,\theta) \in \mathcal{X}_{N}$ with \textit{class} $\mathcal{K}_{\infty}$ function $\hat{\gamma}(\cdot)$. Hence, the optimal rotated objective function is bounded from above and below, according to:
\begin{equation*}
    \gamma(\|x-x^{s}\|) \leq \bar{V}_{N}(x,\theta) \leq \hat{\gamma}(\|x-x^{s}\|).
\end{equation*}
Therefore, $\bar{V}_{N}(\theta_{t},x_{t})$ is a Lyapunov function and $x^{s}$ is an asymptotically stable equilibrium point of \eqref{closedloopsystem} with a region of attraction $\mathcal{X}_{N}$. \qed
\end{pf}

\begin{rem}
We remark that the controlled dissipativity with non parameter-dependent storage functions can also ensure closed-loop convergence and even stability property when the standard dissipativity is not satisfied. The proofs follow similar lines as those of Theorem \ref{thm:convergence} and Corollary \ref{theoreminequality} by fixing the storage function parameter $\theta$ along the system trajectory. Hence, the main benefit of using controlled dissipativity lies in the satisfaction of asymptotic stability, which is unattainable by average constraints alone as pursued in \cite{angeli2011average,muller2014convergence}. Additionally, it is known that verifying the (controlled) dissipativity for any given storage function is a challenging task.Therefore, parameterising the storage function and letting the optimisation solver determine the its values on-line provides a potential simpler approach as used in the EMPC problem \eqref{Opt:main}. Furthermore, exploring more candidates of storage functions within certain classes of function may further improve the economic cost during transient period without loss of the asymptotic stability. In order to verify this conclusion, Section \ref{example:regionalconstraint} will study the benefits of optimising storage function parameters.
\end{rem}

\subsection{Asymptotic performance evaluation}

For systems with guaranteed convergence to the optimal regime of operation, the upper bound of average cost for asymptotically stable closed-loop systems can be obtained easily. In the following, we analyse the asymptotic average performance which is stated as:

\begin{thm}\label{asymptoticperformance}
Suppose Assumption \ref{assum:compactTheta}, \ref{assum:continuousstoragefunction}, \ref{assum:terminal}, and \ref{assumptiondissipative} are satisfied. The asymptotic average performance for closed-loop system using the EMPC controller \eqref{Opt:main} fulfils
\begin{equation}\label{performancebound}
    \limsup_{T \to +\infty} \frac{\sum^{T-1}_{t=0} \ell(x^{*}_{t|t},u^{*}_{t|t}) }{T} \leq \ell(x^{s},u^{s}).
\end{equation}
\end{thm}
\begin{pf}
The proof of this theorem follows similar lines as \cite[Theorem 1]{angeli2011average}. Given the suboptimal and feasible solution at time $t+1$ constructed in the same fashion of \eqref{initialsequence}, the optimal cost-to-go functions fulfil $V_{N}(\theta_{t+1},x_{t+1}) \leq J_{N}(\tilde{x}_{t+1:t+N+1|t+1},\tilde{u}_{t+1:t+N|t+1})$ due to the optimality of $V_{N}(\theta_{t+1},x_{t+1})$, therefore, the following inequality is satisfied
\begin{equation*}
\begin{aligned}
    V_{N}(\theta_{t+1},x_{t+1}) - V_{N}(\theta_{t},x_{t})&\leq  J_{N}(\tilde{x}_{t+1:t+N+1|t+1},\tilde{u}_{t+1:t+N|t+1}) - V_{N}(\theta_{t},x_{t}) \\& = -\ell(x^{*}_{t|t},u^{*}_{t|t}) + V_{f}(\theta^{*}_{t+N+1|t},x^{*}_{t+N+1|t})  - V_{f}(\theta^{*}_{t+N|t},x^{*}_{t+N|t}) + \ell(x^{*}_{t+N|t},\kappa_{f}(x^{*}_{t+N|t}) \\
    & \leq -\ell(x^{*}_{t|t},u^{*}_{t|t}) + \ell(x^{s},u^{s})
\end{aligned}
\end{equation*}
where the inequality holds because of the third condition in Assumption \ref{assum:terminal}.

Then, integrating both sides of the above inequality over any time $T-1 \in \mathbb{I}_{\geq 0}$, dividing by $T$, applying $\liminf$ on both sides, and exploiting
boundedness of solutions, we see that
\begin{equation*}
    0 = \liminf_{T \to \infty} \frac{V_{N}(x_{T})-V_{N}(x_{0})}{T} 
     \leq \liminf_{T \to \infty} \frac{\sum^{T-1}_{t=0} \ell(x^{s},u^{s})- \ell(x^{*}_{t|t},u^{*}_{t|t})}{T} = \ell(x^{s},u^{s}) -  \limsup_{T \to +\infty} \frac{\sum^{T-1}_{t=0} \ell(x^{*}_{t|t},u^{*}_{t|t})}{T}
\end{equation*}
which in turn implies the inequality \eqref{performancebound} and proves the claim in Theorem \ref{asymptoticperformance}. \qed
\end{pf}

\begin{rem}\label{rem:performance}
For more general systems in which the convergence condition, e.g., Assumption \ref{assumptiondissipative}, is not fulfilled, recursive feasibility is lost because the terminal control policy $\kappa_{f}(x)$ might not satisfy the dissipation inequality. However, there always exists a feasible solution due to Assumption \ref{assumptioncontrolleddissipative} and it might happen that the closed-loop system does not converge to the equilibrium. In this case, the upper bound of asymptotic average cost in Theorem \ref{asymptoticperformance} still holds, which will be further explained by an example in Section \ref{sec:example}.
\end{rem}

\section{Examples}\label{sec:example}

In this section, we use numerical examples to visualize the theoretical results. The main goal of this section is to investigate and highlight the validity of a priori convergence to the optimal equilibrium and average cost performance bound when the standard dissipativity property w.r.t the optimal equilibrium is absent.

We consider the following discrete time linear system
\begin{equation*}
    \begin{bmatrix}
    x_{1}\\
    x_{2}
    \end{bmatrix}^{+} = \begin{bmatrix}
    0 & 1\\
    -1 & 0
    \end{bmatrix}
    \begin{bmatrix}
    x_{1}\\
    x_{2}
    \end{bmatrix} + \begin{bmatrix}
    1 \\
    0
    \end{bmatrix} u
\end{equation*}
subject to the state and input coupled constraints
\begin{equation*}
\mathbb{X} = [-1, 1] \times [-1, 1],\quad \mathbb{U}(x) = [-1-x_{2}, 1-x_{2}].
\end{equation*}
The economic stage cost is specified by
\begin{equation*}
    \ell(x,u) = u^{2} + x^{4}_{1} - 0.5x^{2}_{1}.
\end{equation*}
Notice that this cost is not positive definite. In particular, the optimal average performance can be expected to be negative, as the equilibrium solution $x^{s} = [0, 0]^{T}$ is feasible with input $u^{s} = 0$, yielding average cost $\ell(x^{s},u^{s}) = 0$. However, the stage cost can be made negative for some values of $x_{1} \neq 0$.

The zero-input response of the system are (feasible) period $4$ oscillations. The periodic solution $[0.5, 0]^{T}$, $[0, -0.5]^{T}$, $[-0.5, 0]^{T}$, $[0, 0.5]^{T}$ has an average cost $-0.0625$, outperforming the optimal equilibrium $(x^{s},u^{s})$. Therefore, the standard dissipativity property w.r.t the equilibrium does not hold.

Next, we will show the closed-loop system converges to the equilibrium under the proposed EMPC controllers. Thought the following simulations, the prediction horizon is chosen as $N = 20$, the initial system state is $x = [1, 1]^{T}$ and the simulation steps is $T = 100$.

\subsection{Convergence speed and asymptotic average cost}

In the EMPC optimisation problem \eqref{Opt:main}, a terminal equality constraint $x_{t+N|t} = [0, 0]^{T}$ is imposed at every sampling time $t$. The positive definite term for strict controlled dissipativity in \eqref{Opt:dissipativity} is expressed as $\rho(x-x^{s}) = \varrho \cdot \|x-x^{s}\|^{2}_{2}$ where $\varrho \in \mathbb{R}_{\geq 0}$ is a tuning parameter. In addition, the parameter varying storage function is defined in the form of
\begin{equation}\label{storagefunction_TermEqu}
    \lambda(\theta,x) = a_{1}x_{1}^{2} + a_{2}x_{2}^{2} + a_{3}x_{1}x_{2} + a_{4}x_{1} + a_{5}x_{2} + a_{6}
\end{equation}
where $\theta = [a_{1}, a_{2}, a_{3}, a_{4}, a_{5}, a_{6}]^{T} \in \Theta$ with $a_{i} \in [-\bar{\theta}, \bar{\theta}], \forall i \in \mathbb{I}_{[1,6]}$ and
$\Theta = [-\bar{\theta}, \bar{\theta}]^{6}$.

\begin{figure}[H]
\vspace{-3mm}
\centerline{\includegraphics[width=0.6\textwidth]{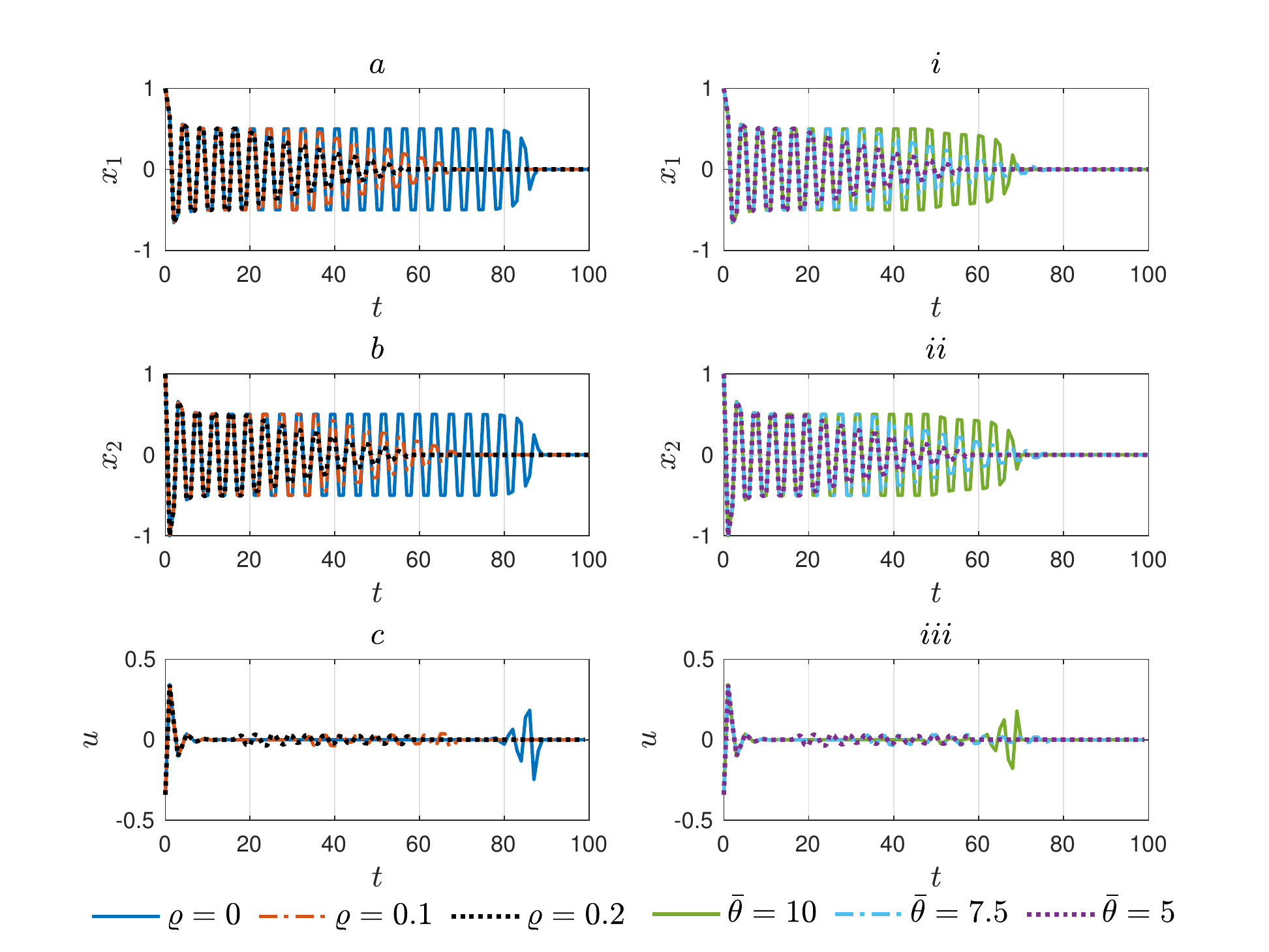}}
\vspace{-5mm}
\caption{Closed-loop state and input trajectories under the EMPC controller \eqref{Opt:main} with different positive definite terms (a, b, c) and different parameter bounds (i, ii, iii) in the dissipation inequality.}
\label{ConvSpeed_TermEqu}
\end{figure}
\begin{figure}[H]
\vspace{-3mm}
\centerline{\includegraphics[width=0.6\textwidth]{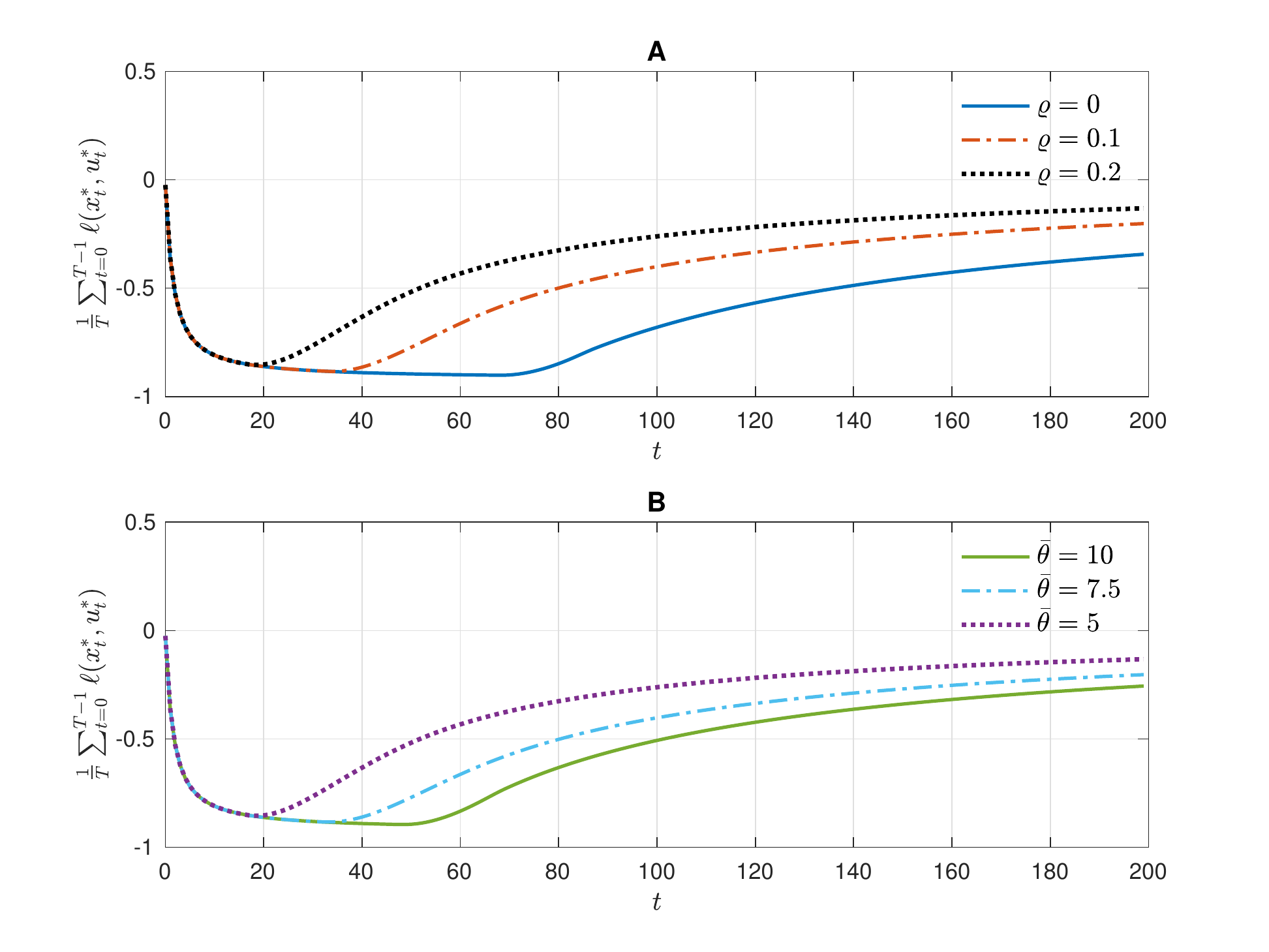}}
\vspace{-5mm}
\caption{Asymptotic average performance under the EMPC controller \eqref{Opt:main} with different A: positive definite term $\varrho$; B: storage function parameter bound $\bar{\theta}$.}
\label{AsyAvCost_TermEqu}
\end{figure}
\begin{figure}[H]
\vspace{-3mm}
\centerline{\includegraphics[width=0.6\textwidth]{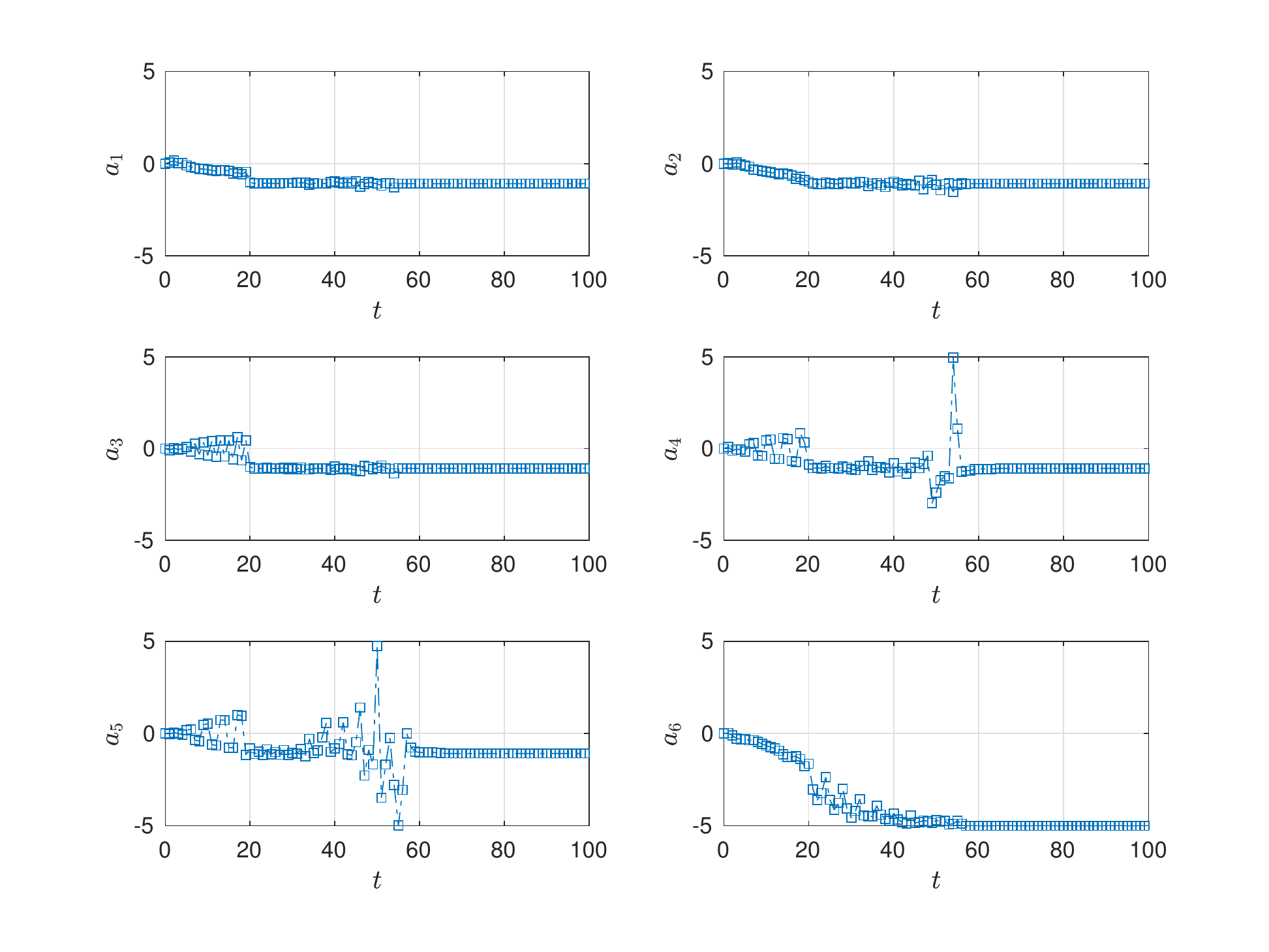}}
\vspace{-5mm}
\caption{Storage function parameters for $\varrho=0.2$ and $\bar{\theta}=5$.}
\label{Theta_TermEqu}
\end{figure}

Under the EMPC controller \eqref{Opt:main}, the state and input components converge to the optimal steady-state asymptotically. Fig.~\ref{ConvSpeed_TermEqu} compares the achieved results of varying parameter bounds and weights in the positive definite term. Specifically, when the parameter upper bound is fixed at $\bar{\theta} = 5$, the results of tuning $\varrho$ are shown in the subplots Fig.~\ref{ConvSpeed_TermEqu} (a,b,c) to see the impact of adding the positive definite term on convergence speed. As $\varrho$ is increased from $0$ to $0.2$, it can be seen that the state trajectory starts dampening the periodic behavior earlier and converges to the optimal equilibrium faster. As a result, a loss in transient average profit, namely, $\frac{1}{T}\sum^{T}_{t=0}\ell(x^{*}_{t},u^{*}_{t})$, is observed by increasing $\varrho$ as shown in Fig.~\ref{AsyAvCost_TermEqu}.A. On the other hand, given $\varrho = 0.2$, faster convergence of the state trajectory is obtained by decreasing $\bar{\theta}$ as observed from the subplots Fig.~\ref{ConvSpeed_TermEqu} (i,ii,iii). However, according to Fig.~\ref{AsyAvCost_TermEqu}.B, the tighter bound of storage function parameters will also lose more economic performance over the transient period. Therefore, we can prove that the steady-state operation is suboptimal in the sense of economic profit and, moreover, conclude a trade-off between economic cost and convergence performance due to the effect of parameters $\varrho$ and $\bar{\theta}$. Fig.~\ref{Theta_TermEqu} shows the profiles of parameters in the polynomial storage function \eqref{storagefunction_TermEqu}. In particular, the constant term $a_{6}$ is decreasing over the transient phase until hitting the lower bound to provide feasibility for the dissipation inequality and strive for more economic profit. When the storage function parameters cannot be further decreased, the system states have to converge to the equilibrium for suboptimal operation as a result of the controlled dissipativity.

Regarding the argument in Remark \ref{rem:performance}, it can be verified by considering the same system state equation as the above but a different stage cost function, e.g., $\ell(x,u) = |x_{1}x_{2}|$. For this stage cost, there is one equilibrium $x^{s} = [0, 0]^{T}$ and infinite number of rotating solutions $x = [c, 0]^{T}$ or $x = [0, c]^{T}$ where $c$ is an arbitrary constant within the interval $[-1, 1]$. For instance, let the system state is initiated at $[x_{1}, 0]^{T}$ where $x_{1}>0$, the closed-loop system may keep rotating rather than converging to the equilibrium $x^{s}$ as the supply rate is always $s(x,u)=0$, and the average cost is equal to $\ell(x^{s},u^{s})$ which validates the asymptotic performance upper bound \eqref{performancebound} in Theorem \ref{asymptoticperformance}.

\subsection{Imposing stability condition}\label{example:regionalconstraint}

In the following simulation, we choose the terminal penalty function by $V_{f}(x) = x_{2}^{2}, \, \forall x \in \mathbb{X}_{f}$ with the terminal region $\mathbb{X}_{f} = \{x \in \mathbb{X}\,|\, x_{1} = 0, x_{2} \in [-1, 1]\}$. The local control policy is selected as $\kappa_{f}(x) = -x_{2}$ such that the conditions in Assumption \ref{assum:terminal} are fulfilled. Meanwhile, for all state $x \in \mathbb{X}_{f}$, this local policy also ensures the satisfaction of Assumption \ref{assumptiondissipative}, viz,
\begin{equation*}
    \lambda(\theta^{+},f(x,\kappa_{f}(x))) - \lambda(\theta,x) \leq \ell(x,u) - \ell(x^{s},u^{s}) = x_{2}^{2}.
\end{equation*}
Then, the storage function is selected in the polynomial form of
\begin{equation*}
    \lambda(\theta,x) = a_{1}x_{1}^{4} + a_{2}x_{1}^{3} + a_{3}x_{1}^{2} + a_{4}x_{1} + a_{5}
\end{equation*}
which lead to $\lambda(\theta,x) = a_{5}, \forall x \in \mathbb{X}_{f}$. A candidate of $\theta^{+}$ at the tail of the prediction horizon for ensuring recursive feasibility is $\theta^{+} = \theta$. In order to satisfy Assumption \ref{assum:stability}, it requires $a_{5} = 0$ such that $\lambda(\theta,x^{s}) = 0, \forall a_{i} \in [-5, 5]$ where $i \in \mathbb{I}_{[1,4]}$. The positive definite term in the dissipation inequality is expressed as $\rho(x-x^{s}) = 0.2 \cdot \|x-x^{s}\|^{2}_{2}$.

\begin{figure}[H]
\vspace{-3mm}
\centerline{\includegraphics[width=0.6\textwidth]{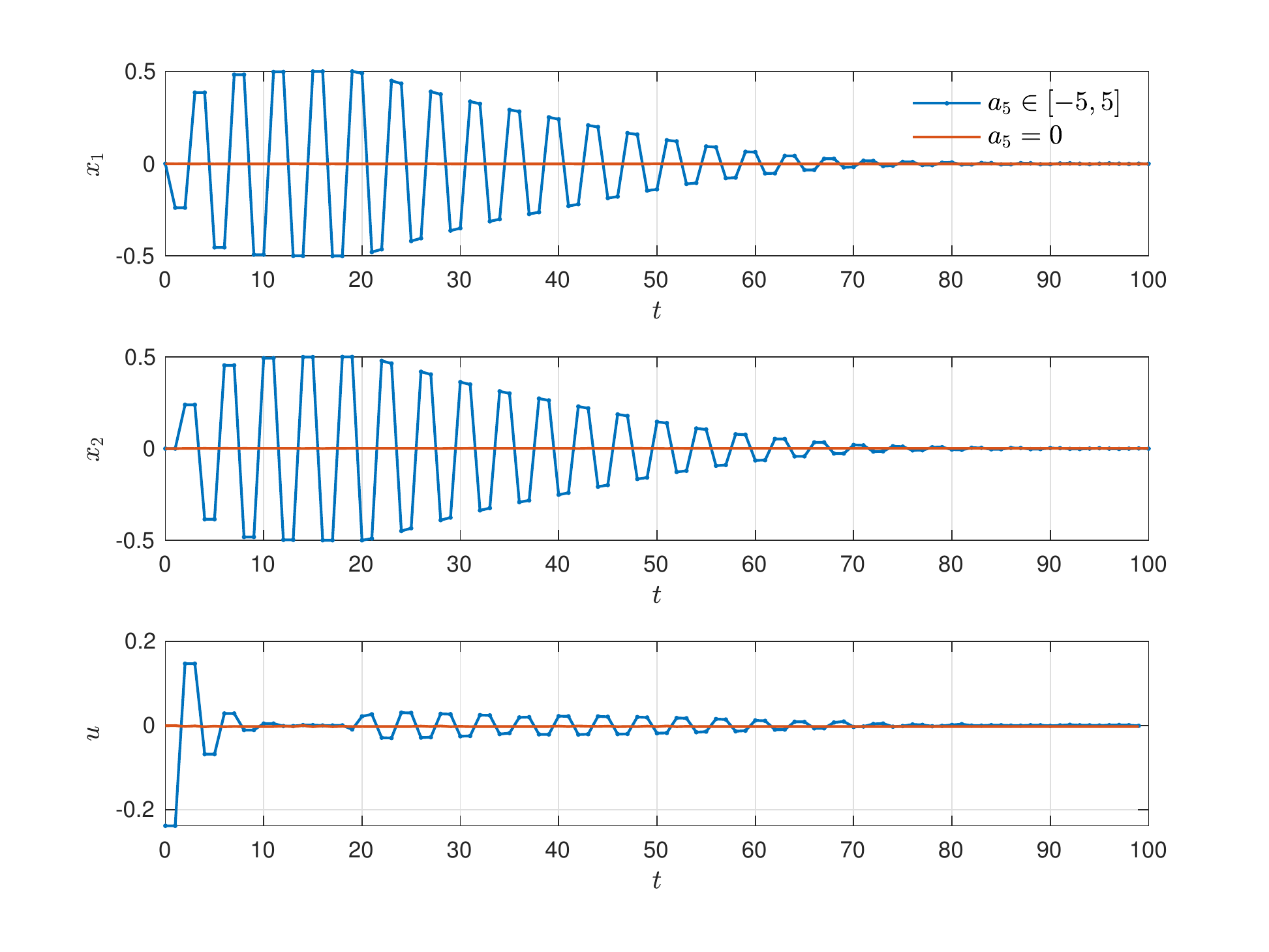}}
\vspace{-5mm}
\caption{Closed-loop system behavior with/without the stability condition in Assumption \ref{assum:stability}.}
\label{Convergence_vs_Stability}
\end{figure}

The simulation results to show effectiveness of Assumption \ref{assum:stability} on closed-loop stability are presented in Fig.~\ref{Convergence_vs_Stability}. When the condition in Assumption \ref{assum:stability} is absent, i.e., $a_{5} \in [-5, 5]$, the state trajectory starting from the optimal steady-state $x^{s}$ at $t=0$ exhibits oscillations because $x^{s}$ is not stable and an improved (transient) performance of the closed-loop system can be achieved moving away from equilibrium. When Assumption \ref{assum:stability} is imposed, i.e., $a_{5} = 0$, the closed-loop system stays at $x^{s}$ which means that $x^s$ is a stable equilibrium point of the closed-loop system.

\begin{figure}[H]
\vspace{-3mm}
\centerline{\includegraphics[width=0.6\textwidth]{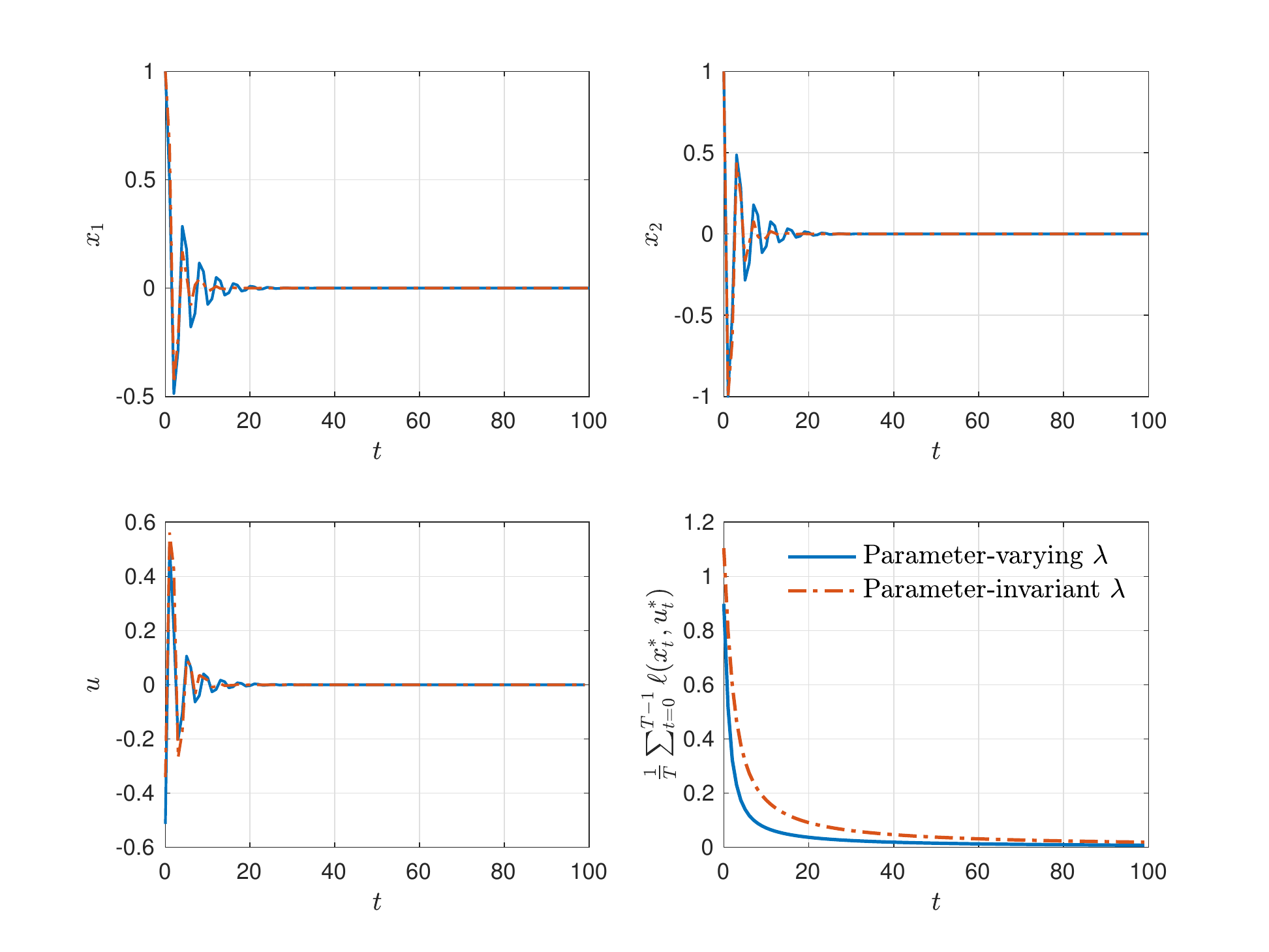}}
\vspace{-5mm}
\caption{Closed-loop system trajectories and the asymptotic average cost when the parameter-varying and parameter-invariant storage functions are used.}
\label{Varying_vs_Invariant}
\end{figure}

In order to verify the economic benefits of setting storage function parameters as decision variables, two simulations with the same initial state $x_{0} = [1, 1]^{T}$ are implemented. In the first simulation, the storage function is selected as $\lambda(\theta,x) = 0, \forall x \in \mathbb{X}$, i.e., $a_{i}=0, \forall i \in \mathbb{I}_{[1,5]}$, which satisfies the controlled dissipativity as the supply rate $\ell(x,u)-\ell(x^{s},u^{s}) = x^{2}_{2}$ is always non-negative. In the second simulation, we choose $a_{i} \in [-5, 5], \forall i \in \mathbb{I}_{[1,4]}$ and $a_{5} = 0$. Based on the developed theoretical analysis, both selections of the storage function will ensure the closed-loop system asymptotic stability with respect to the optimal steady-state, which is also demonstrated by the resulting trajectories in Fig.~\ref{Varying_vs_Invariant}. Concerning the asymptotic average cost in Fig.~\ref{Varying_vs_Invariant}, it is observed that using the storage function parameters as optimisation variables can achieve lower economic cost over the transient period. Hence, we conclude that the controlled dissipativity under certain storage function condition can ensures stability property, and moreover, the adoption of storage function parameters as optimisation variables can achieve improved economic performance over the transient period.


\section{Conclusion}\label{sec:conclusion}

This paper presented a novel concept of controlled dissipativity with parameter varying storage functions. Based on the controlled dissipation inequality, EMPC controllers are designed to drive the system to a steady-state without requiring the standard dissipativity condition w.r.t it. Recursive feasibility is ensured by iteratively constructing warm-start for the optimisation problems. Theoretical analysis of the asymptotic stability for the best equilibrium and the asymptotic average performance are addressed showing the main contribution of this current work. The practical application of the presented approaches were illustrated by means of numerical examples indicating the trade-off between the convergence speed and the average economic performance. The effectiveness of imposing certain condition on storage function to ensure closed-loop stability is also evaluated.

\bibliographystyle{plain}        
\bibliography{autosam}           

\begin{thebibliography}{10}

\bibitem{amrit2011economic}
Rishi Amrit, James~B Rawlings, and David Angeli.
\newblock Economic optimization using model predictive control with a terminal
  cost.
\newblock {\em Annual Reviews in Control}, 35(2):178--186, 2011.

\bibitem{angeli2011enforcing}
David Angeli, Rishi Amrit, and James~B Rawlings.
\newblock Enforcing convergence in nonlinear economic {MPC}.
\newblock In {\em 2011 50th IEEE Conference on Decision and Control and
  European Control Conference}, pages 3387--3391. IEEE, 2011.

\bibitem{angeli2011average}
David Angeli, Rishi Amrit, and James~B Rawlings.
\newblock On average performance and stability of economic model predictive
  control.
\newblock {\em IEEE transactions on automatic control}, 57(7):1615--1626, 2011.

\bibitem{angeli2016theoretical}
David Angeli, Alessandro Casavola, and Francesco Tedesco.
\newblock Theoretical advances on economic model predictive control with
  time-varying costs.
\newblock {\em Annual Reviews in Control}, 41:218--224, 2016.

\bibitem{bayer2014tube}
Florian~A Bayer, Matthias~A M{\"u}ller, and Frank Allg{\"o}wer.
\newblock Tube-based robust economic model predictive control.
\newblock {\em Journal of Process Control}, 24(8):1237--1246, 2014.

\bibitem{berberich2018indefinite}
Julian Berberich, Johannes K{\"o}hler, Frank Allg{\"o}wer, and Matthias~A
  M{\"u}ller.
\newblock Indefinite linear quadratic optimal control: Strict dissipativity and
  turnpike properties.
\newblock {\em IEEE Control Systems Letters}, 2(3):399--404, 2018.

\bibitem{berberich2020dissipativity}
Julian Berberich, Johannes K{\"o}hler, Frank Allg{\"o}wer, and Matthias~A
  M{\"u}ller.
\newblock Dissipativity properties in constrained optimal control: a
  computational approach.
\newblock {\em Automatica}, 114:108840, 2020.

\bibitem{diehl2010lyapunov}
Moritz Diehl, Rishi Amrit, and James~B Rawlings.
\newblock A {L}yapunov function for economic optimizing model predictive
  control.
\newblock {\em IEEE Transactions on Automatic Control}, 56(3):703--707, 2010.

\bibitem{dong2018analysis}
Zihang Dong and David Angeli.
\newblock Analysis of economic model predictive control with terminal penalty
  functions on generalized optimal regimes of operation.
\newblock {\em International Journal of Robust and Nonlinear Control},
  28(16):4790--4815, 2018.

\bibitem{dong2020homothetic}
Zihang Dong and David Angeli.
\newblock Homothetic tube-based robust economic {MPC} with integrated moving
  horizon estimation.
\newblock {\em IEEE Transactions on Automatic Control}, 66(1):64--75, 2020.

\bibitem{ellis2014tutorial}
Matthew Ellis, Helen Durand, and Panagiotis~D Christofides.
\newblock A tutorial review of economic model predictive control methods.
\newblock {\em Journal of Process Control}, 24(8):1156--1178, 2014.

\bibitem{faulwasser2018economic}
Timm Faulwasser, Lars Gr{\"u}ne, Matthias~A M{\"u}ller, et~al.
\newblock Economic nonlinear model predictive control.
\newblock {\em Foundations and Trends{\textregistered} in Systems and Control},
  5(1):1--98, 2018.

\bibitem{faulwasser2014turnpike}
Timm Faulwasser, Milan Korda, Colin~N Jones, and Dominique Bonvin.
\newblock Turnpike and dissipativity properties in dynamic real-time
  optimization and economic {MPC}.
\newblock In {\em 53rd ieee conference on decision and control}, pages
  2734--2739. IEEE, 2014.

\bibitem{grune2013economic}
Lars Gr{\"u}ne.
\newblock Economic receding horizon control without terminal constraints.
\newblock {\em Automatica}, 49(3):725--734, 2013.

\bibitem{grune2018turnpike}
Lars Gr{\"u}ne and Roberto Guglielmi.
\newblock Turnpike properties and strict dissipativity for discrete time linear
  quadratic optimal control problems.
\newblock {\em SIAM Journal on Control and Optimization}, 56(2):1282--1302,
  2018.

\bibitem{grune2016relation}
Lars Gr{\"u}ne and Matthias~A M{\"u}ller.
\newblock On the relation between strict dissipativity and turnpike properties.
\newblock {\em Systems \& Control Letters}, 90:45--53, 2016.

\bibitem{grune2014asymptotic}
Lars Gr{\"u}ne and Marleen Stieler.
\newblock Asymptotic stability and transient optimality of economic {MPC}
  without terminal conditions.
\newblock {\em Journal of Process Control}, 24(8):1187--1196, 2014.

\bibitem{heidarinejad2012economic}
Mohsen Heidarinejad, Jinfeng Liu, and Panagiotis~D Christofides.
\newblock Economic model predictive control of nonlinear process systems using
  {L}yapunov techniques.
\newblock {\em AIChE Journal}, 58(3):855--870, 2012.

\bibitem{kohler2018periodic}
Johannes K{\"o}hler, Matthias~A M{\"u}ller, and Frank Allg{\"o}wer.
\newblock On periodic dissipativity notions in economic model predictive
  control.
\newblock {\em IEEE Control Systems Letters}, 2(3):501--506, 2018.

\bibitem{lucia2014handling}
Sergio Lucia, Joel~AE Andersson, Heiko Brandt, Moritz Diehl, and Sebastian
  Engell.
\newblock Handling uncertainty in economic nonlinear model predictive control:
  A comparative case study.
\newblock {\em Journal of Process Control}, 24(8):1247--1259, 2014.

\bibitem{marquez2014min}
Alejandro Marquez, Julian Pati{\~n}o, and Jairo Espinosa.
\newblock Min-max economic model predictive control.
\newblock In {\em 53rd IEEE Conference on Decision and Control}, pages
  4410--4415. IEEE, 2014.

\bibitem{muller2014necessity}
Matthias~A M{\"u}ller, David Angeli, and Frank Allg{\"o}wer.
\newblock On necessity and robustness of dissipativity in economic model
  predictive control.
\newblock {\em IEEE Transactions on Automatic Control}, 60(6):1671--1676, 2014.

\bibitem{muller2014convergence}
Matthias~A M{\"u}ller, David Angeli, Frank Allg{\"o}wer, Rishi Amrit, and
  James~B Rawlings.
\newblock Convergence in economic model predictive control with average
  constraints.
\newblock {\em Automatica}, 50(12):3100--3111, 2014.

\bibitem{muller2016economic}
Matthias~A M{\"u}ller and Lars Gr{\"u}ne.
\newblock Economic model predictive control without terminal constraints for
  optimal periodic behavior.
\newblock {\em Automatica}, 70:128--139, 2016.

\bibitem{pirkelmann2019approximate}
Simon Pirkelmann, David Angeli, and Lars Gr{\"u}ne.
\newblock Approximate computation of storage functions for discrete-time
  systems using sum-of-squares techniques.
\newblock {\em IFAC-PapersOnLine}, 52(16):508--513, 2019.

\bibitem{rawlings2012fundamentals}
James~B Rawlings, David Angeli, and Cuyler~N Bates.
\newblock Fundamentals of economic model predictive control.
\newblock In {\em 2012 IEEE 51st IEEE conference on decision and control
  (CDC)}, pages 3851--3861. IEEE, 2012.

\bibitem{rawlings2008unreachable}
James~B Rawlings, Dennis Bonn{\'e}, John~B Jorgensen, Aswin~N Venkat, and
  Sten~Bay Jorgensen.
\newblock Unreachable setpoints in model predictive control.
\newblock {\em IEEE Transactions on Automatic Control}, 53(9):2209--2215, 2008.

\bibitem{rawlings2017model}
James~Blake Rawlings, David~Q Mayne, and Moritz Diehl.
\newblock {\em Model predictive control: theory, computation, and design},
  volume~2.
\newblock Nob Hill Publishing Madison, WI, 2017.

\bibitem{zanon2016periodic}
Mario Zanon, Lars Gr{\"u}ne, and Moritz Diehl.
\newblock Periodic optimal control, dissipativity and {MPC}.
\newblock {\em IEEE Transactions on Automatic Control}, 62(6):2943--2949, 2016.

\end{thebibliography}



\end{document}